\documentclass[a4,12pt]{scrartcl}

\usepackage{xcolor} 
\definecolor{ocre}{RGB}{52,177,201} 

\usepackage{diagbox}
\usepackage{ stmaryrd }

\usepackage{listings}
\usepackage{afterpage}

\usepackage{avant} 
\usepackage{mathptmx} 

\usepackage{microtype} 
\usepackage[utf8]{inputenc} 
\usepackage[T1]{fontenc} 

\usepackage{float}
\usepackage{caption}
\usepackage{subcaption}
\usepackage[final]{graphicx}

\usepackage{tikz}
\usepackage{pgfplots} 
\usepackage{pgfplotstable}
\usetikzlibrary{math}
\usetikzlibrary{calc}
\usepackage{tkz-euclide}
\usetikzlibrary{patterns}
\usetikzlibrary{decorations.pathmorphing,shapes.geometric,arrows}
\usepgfplotslibrary{external} 

\usepackage{mdwlist}
\usepackage{multirow}
\usepackage{siunitx}

\usepackage[edges]{forest}

\def\Size{4pt}
\tikzset{
	folder/.pic={
		\filldraw[draw=folderborder,top color=folderbg!50,bottom color=folderbg]
		(-1.05*\Size,0.2\Size+5pt) rectangle ++(.75*\Size,-0.2\Size-5pt);  
		\filldraw[draw=folderborder,top color=folderbg!50,bottom color=folderbg]
		(-1.15*\Size,-\Size) rectangle (1.15*\Size,\Size);
	}
}

\usepackage{graphicx} 
\graphicspath{{img/}} 

\usepackage{lipsum} 

\usepackage{tikz} 

\usepackage[english]{babel} 

\usepackage{enumitem} 
\setlist{nolistsep} 

\usepackage{booktabs} 

\usepackage{eso-pic} 

\usepackage[titletoc]{appendix} 

\definecolor{folderbg}{RGB}{124,166,198}
\definecolor{folderborder}{RGB}{110,144,169}
\definecolor{foldercolor}{RGB}{124,166,198}

\definecolor{folderbg}{RGB}{124,166,198}
\definecolor{folderborder}{RGB}{110,144,169}

\definecolor{Navy}{RGB}{000,000,128}
\definecolor{DarkRed}{RGB}{139,000,000}
\definecolor{DarkOrange}{RGB}{205,102,000}
\definecolor{DarkGreen}{RGB}{000,100,000}
\definecolor{Beige}{RGB}{250,250,240}
\definecolor{applegreen}{rgb}{0.55, 0.71, 0.0}
\definecolor{palecopper}{rgb}{0.85, 0.54, 0.4}
\definecolor{platinum}{rgb}{0.9, 0.89, 0.89}
\definecolor{CAUpurple}{RGB}{153,000,153}
\definecolor{CAUblue}{RGB}{051,051,153}
\definecolor{DarkPurple}{rgb}{0.07, 0.04, 0.56}

\definecolor{U0}{RGB}{183,028,028}
\definecolor{U1}{RGB}{136,014,079}
\definecolor{U2}{RGB}{074,020,140}
\definecolor{U3}{RGB}{049,027,146}
\definecolor{U4}{RGB}{026,035,126}
\definecolor{U5}{RGB}{013,071,161}
\definecolor{U6}{RGB}{001,087,155}
\definecolor{U7}{RGB}{000,096,100}
\definecolor{U8}{RGB}{000,077,064}
\definecolor{U9}{RGB}{027,094,032}

\definecolor{codegreen}{rgb}{0,0.6,0}
\definecolor{codegray}{rgb}{0.5,0.5,0.5}
\definecolor{codepurple}{rgb}{0.58,0,0.82}
\definecolor{backcolour}{rgb}{0.95,0.95,0.92}

\lstdefinestyle{mystyle}{
	backgroundcolor=\color{backcolour},   
	commentstyle=\color{codegreen},
	keywordstyle=\color{magenta},
	numberstyle=\tiny\color{codegray},
	stringstyle=\color{codepurple},
	basicstyle=\ttfamily\footnotesize,
	breakatwhitespace=false,         
	breaklines=true,                 
	captionpos=b,                    
	keepspaces=true,                 
	numbers=left,                    
	numbersep=5pt,                  
	showspaces=false,                
	showstringspaces=false,
	showtabs=false,                  
	tabsize=2
}

\usepackage{fancyhdr} 

\usepackage{amsmath,amsfonts,amssymb,amsthm} 

\usepackage[backend=bibtex,
citestyle=numeric]{biblatex}

\bibliography{sources/sources}

\usepackage{hyperref}

\hypersetup{
    bookmarks=true,
    unicode=false,
    pdftoolbar=true,
    pdfmenubar=true,
    pdffitwindow=false,
    pdfstartview={FitH},
    pdftitle={My title},
    pdfauthor={B. Philippi},
    pdfsubject={Subject},
    pdfcreator={Creator},
    pdfproducer={Producer},
    pdfkeywords={keyword1} {key2} {key3},
    pdfnewwindow=true,
    colorlinks=false,
    linkbordercolor={0 0 1},
    linkcolor=blue,
    citecolor=green,
    filecolor=magenta,
    urlcolor=cyan
}

\theoremstyle{definition}

\usepackage[left=2cm,right=2cm,top=2cm,bottom=2cm,includeheadfoot]{geometry}

\title{A Micro-Macro parallel-in-time Implementation for the 2D Navier-Stokes Equations}
\author{Benedict $\text{Philippi}^1$, Mahfuz $\text{Sarker Miraz}^2$ , Thomas $\text{Slawig}^3$} 
\date{\empty}

\usepackage{csquotes}

\begin{document}

\maketitle

\begin{abstract}
	In this paper the Micro-Macro Parareal algorithm as proposed in \cite{Samaey2012} was adapted to PDEs. The parallel-in-time approach requires two meshes of different spatial resolution in order to compute approximations in an iterative way to a predefined reference solution. When fast convergence in few iterations can be accomplished the algorithm is able to generate wall-time reduction in comparison to the serial computation. We chose the laminar flow around a cylinder benchmark on 2-dimensional domain which was simulated with the open-source software OpenFoam. The numerical experiments presented in this work aim to approximate states local in time and space and the diagnostic lift coefficient. The Reynolds number is gradually increased from 100 to 1,000, before the transition to turbulent flows sets in. After the results are presented the convergence behavior is discussed with respect to the Reynolds number and the applied interpolation schemes.
\end{abstract}

\vspace{0.25cm}
\small
\textbf{${}^1$ \textit{Christian-Albrecht-Universität Kiel, Dept. of Computer Science, b.k.philippi@gmail.com}}

\textbf{${}^2$ \textit{Christian-Albrecht-Universität Kiel, Dept. of Computer Science}}

\textbf{${}^3$ \textit{Christian-Albrecht-Universität Kiel, Dept. of Computer Science, ts@informatik.uni-kiel.de}}
\normalsize

\newpage

\section{Introduction}

The application of spatial domain decomposition techniques to reduce wall-times for computational fluid dynamics (CFD) allows for significant reduction in run-times. With more cores than ever being available on modern high performance computers (HPC) the speed-up does not grow accordingly. Run-time reduction due to domain decomposition saturates at some point, even though far more processors could be included in the process. The application of parallel-in-time algorithms aims to overcome this saturation point and generate additional speed-ups by the concept of time-parallelism. In the field of CFD the high demand for spatial and temporal resolution is computationally expensive and even with domain decomposition techniques further speed-ups are desirable. One possibility for parallelization in time is the Parareal algorithm, introduced by \cite{Maday2001} in 2001. Since then it has been applied to a broad variety of problems, compare \url{parallel-in-time.org}. Its non-intrusive nature makes it suitable for application to all kinds of initial value problems (IVP). In this study we aim for the time-parallel solution of the 2D Navier-Stokes equations. We adapted the Micro-Macro Parareal algorithm by \cite{Samaey2012} to the laminar flow around a cylinder test case for a variety of Reynolds numbers. Attempts parallelize the cylinder flow in time had been made by \cite{KreienbuehlEtAl2015} for different numerical time integration schemes and with a two mesh (Micro-Macro) approach in \cite{FischerEtAl2005}. Other studies towards Parareal convergence for the Navier-Stokes equations are given by \cite{SteinerEtAl2015} for the 2D driven cavity test case and in 3D by \cite{Ruprecht2014}. The results of those studies show the impact of spatial resolution and the Reynolds number on the convergence behavior, as Parareal starts struggling with the increase of both. On the other hand it could be shown that for the time-parallel approximation of diagnostic variables the concept of coarsening in space can be successfully applied in the case of homogeneous turbulence, see \cite{Lunet2018}. 

In this study we aim for the diagnostic lift coefficient of the laminar flow around the cylinder. We investigate the needed convergence in prognostic variables, like velocity, in order to compute a reasonable time history of the lift coefficient. The Micro-Macro Parareal algorithm was originally designed for ODE initial value problems and later implemented for an Energy Balance Model (EBM) in 1D \cite{Slawig2018}. The coarse propagator was defined by a scalar macro-scale ODE (0-D EBM) and the fine propagator for the temperature distribution in latitudinal direction in 1D on the northern sphere. The approach suggested a self-consistent matching operator for the micro- and macro-scale. Further, demands were made on the lifting and restriction interpolation schemes used. As the matching operator was introduced for a different setup in \cite{Samaey2012}\cite{Slawig2018}, we have to adapt it to the present setting in compliance with the proposed consistency requirements. In this work, the focus was placed on the feasibility of the time-parallel approach and thus, all simulations are carried out on a single core, leaving the concept of full space-time-parallelism for future investigations. We will discuss the results with respect to the used interpolation schemes and the algorithm's ability to converge to the lift coefficient in dependency of the Reynolds number.

\section{The parallel-in-time Algorithm}

\subsection{Parareal}

The Parareal algorithm aims to decompose the time domain and solve the IVP on each part in parallel in an iterative way. Therefore, a given simulation time interval $[0,T]$ is decomposed into time slices of equal size:
\begin{equation}
	0 \; = \; t_0 \; < \; t_1 \; < \; \dots \; < \; t_N \; = T \; .
\end{equation}
Each time slice is defined by $\Delta T_n = [t_n,t_{n+1}]$ with $n=0,\dots,N_t-1$. With Parareal being an iterative solver, we introduce two propagators of different computational cost. A fast solver $G$ and a fine solver $F$. The latter is considered to be computationally expensive and represents the reference for solving the problem. Both propagators are solving the IVP defined on each time slice $\Delta T_n$ with an initial value $U_n$. Throughout this paper the following notation will be used:
\begin{equation}
	\begin{split}
		G_{U_n,t_n,t_{n+1}} \; &= \; G_{n+1}(U_n) \; , \\ 
		F_{U_n,t_n,t_{n+1}} \; &= \; F_{n+1}(U_n) \; .
	\end{split}
\end{equation}
The propagators can differ in the numerical time integration method and step size. The step sizes for the coarse and fine solver are defined as $\Delta t_C$ and $\Delta t_F$, respectively. Parareal aims to execute the fine solver $F$ on the time slices $\Delta T_n$ in parallel and therefore requires initial values for each. The computationally cheap propagator $G$ will provide these values in a serial run for all time slices. The execution of $F$ and $G$ in an iterative way is defined by the Parareal algorithm as:
\begin{equation}
	\begin{split}
		U^{k+1}_{n+1} \; &= \; G^{k+1}(U^{k+1}_n) \; + \; F^k_{n+1}(U^k_n) \; - \; G^k_{n+1}(U^k_n) \; , \quad k = 0,\dots,N_t-1 \; ,\\
		U^0_{n+1} \; &= \; G^0(U^0_n) \; , \\
		U^k_0 \; &= \; u_0 .
	\end{split}
\end{equation}
Here, $k$ denotes the iteration index. The 0th iteration is defined by the coarse propagator $G$ and the initial condition $u_0$ is given by the problem to be solved. The iterative solution $U^{k+1}_{n+1}$ will converge eventually to the serial fine approximation at $t_{n+1}$:
\begin{equation}
	U^{k+1}_{n+1} \; = \; F_{n+1} \; .
\end{equation}
The iteration will necessarily converge to the fine reference solution for $k=N_t-1$. In order to generate speed-ups Parareal is required to meet a pre-defined error threshold $\varepsilon$ in $k \ll N_t-1$ iterations.

\subsection{Micro-Macro Parareal}

The Micro-Macro Parareal algorithm was introduced in \cite{Samaey2012} for multi-scale ODEs. The concept considered a coarse propagator that solves the problem on a slow macroscopic scale. In a time-parallel iteration the coarse solution is updated by a microscopic solver to approximate the chosen multi-scale ODE. In \cite{Slawig2018} the approach was applied to an Energy-Balance-Model (EBM) in 1-D that computes approximations to the temperature on the upper half of a sphere with the assumption of symmetry at the equator and in longitude. A simple 0-D model was chosen as the coarse propagator to provide a fast solver for serial computations. The Parareal approximations are defined on two different spatial domains. Thus, the update procedure on each time slice is given as a two-step method. The iterative solution by Parareal on the macroscopic scale is defined as $\hat{U}$ and $U$ for the target microscopic scale. The Micro-Macro Parareal iteration reads:
\begin{equation}
	\begin{split}
		\hat{U}^{k+1}_{n+1} \; &= \; \hat{G}^{k+1}_{n+1} \; + \; \mathbf{R}(F^k_{n+1}) \; - \; \hat{G}^k_{n+1} \; , \\
		U^{k+1}_{n+1} \; &= \; \mathbf{P} (\hat{U}^{k+1}_{n+1},F^k_{n+1}) \; .
	\end{split}
\end{equation}
For the sake of clarity the terms have been further abbreviated. The algorithm requires two operators to map the propagators solution to the respective other domain. The restriction operator $\mathbf{R}$ for interpolation from fine to coarse and a lifting operator $\mathbf{L}$ for vice versa. The approach by \cite{Slawig2018} utilizes a matching operator $\mathbf{P}$, which matches a macroscopic with a microscopic solution. In contrast to $\mathbf{L}$, the matching operator $\mathbf{P}$ needs microscopic information in addition to the coarse state as an input. The restriction operator $\mathbf{R}$ computes an average of the fine 1-D solution for use in the 0-D coarse solver. Then, the operator $\mathbf{P}$ was defined by:
\begin{equation}
	\mathbf{P}(\hat{U}^{k+1}_{n+1},F^k_{n+1}) \; = \; \hat{U}^{k+1}_{n+1} \frac{F^k_{n+1}}{\mathbf{R}(F^k_{n+1})} \; , \quad \text{with} \quad F^k_{n+1} \; \in \mathbb{R}^d  \; , \; \hat{U}^{k+1}_{n+1} \; , \; \mathbf{R}(F^k_{n+1}) \; \in \mathbb{R} \; .
	\label{EQ:P}
\end{equation}
The lift operator $\mathbf{L}$ is applied after the 0-th iteration when no fine scale information is available yet. The interpolation operators have to fulfill the following conditions:
\begin{equation}
	\begin{split}
		1.& \quad \mathbf{R} \circ \mathbf{L} (\hat{U}) \; = \; Id. \\
		2.& \quad (\mathbf{R} \circ \mathbf{P})(\hat{U},F) \; = \; \hat{U} \\
		3.& \quad \mathbf{P}(\mathbf{R}(F),F) \; = \; F 
	\end{split}
	\label{EQ:requirements}
\end{equation}
The second and third requirement are consistency properties of the matching operator $\mathbf{P}$, compare \cite{Samaey2012}. The second requirement represents the equivalent of the desired property for the lifting operator $\mathbf{L}$ in the first requirement. 

Since both propagators in this work are defined by the same PDEs on the same spatial domain with different resolution, the algorithm must be adapted. The matching operator $\mathbf{P}$ in Eq.\ref{EQ:P} is not applicable in a meaningful way when both solutions have vector-valued states instead of a scalar quantity, like temperature. In order to adapt the operator while complying with the demand for consistency in Eq.\ref{EQ:requirements} the Full Approximation Scheme (FAS) for non-linear problems from the multigrid methods \cite{FAS} was taken up:
\begin{equation}
	\begin{split}
		\hat{U}^{k+1}_{n+1} \; &= \; \hat{G}^{k+1}_{n+1} \; + \; \mathbf{R}(F^k_{n+1}) \; - \; \hat{G}^k_{n+1} \; , \\
		U^{k+1}_{n+1} \; &= \; \mathbf{L}(\hat{U}^{k+1}_{n+1}) \; + \; F^k_{n+1} \; - \; \mathbf{L} \circ \mathbf{R} (F^k_{n+1}) \; .
	\end{split}
	\label{EQ:MMPAR}
\end{equation}
To account for the loss of information during the update procedure on the macroscopic scale a micro-scale contribution is added by the expression $F^k_{n+1} - \mathbf{L} \circ \mathbf{R} (F^k_{n+1})$ to the lifted iterative solution $\mathbf{L}(\hat{U}^{k+1}_{n+1})$. Further, we are able to preserve the consistency properties in Eq.\ref{EQ:requirements}. The requirements are checked for different interpolation methods in Tabs.\ref{TAB:intpRe100} and \ref{TAB:intpRe1000}, after the test case has been introduced. 

\subsection{A priori Speedup Estimate}

The introduction of the Micro-Macro Parareal algorithm concludes with an a priori speed-up estimate. Although a reduction of runtimes should always be done by wall-time measurements, a prior estimation of the speedup can be utilized to determine an optimal operational regime of Parareal. The speed-up estimate is given by:
\begin{equation}
	S \; = \; \min \left(\frac{m}{k+1} , \frac{N_t}{k} \right) \; , \quad \text{with} \quad m \; = \; \frac{\tau_F}{\tau_C} \; ,
	\label{EQ:speedup}
\end{equation}
where $m$ denotes the run-time ratio of the fine and coarse propagator $\tau_F$ and $\tau_C$ for are given time interval. In order to generate speed-ups the algorithm is required to converge in few iterations, such that $k \ll m$ and $k \ll N_t$. In Eq.\ref{EQ:speedup} the wall-time necessary for interpolation and the iteration procedure is not considered. Thus, the estimate can be seen as an upper bound on the expected speed-up with respect to the amount of iterations. We used Eq.\ref{EQ:speedup} to limit the amount of iterations to a threshold $K$ for the numerical experiments in Sec.\ref{SEC:EXP}.

\section{Test Case}

During this work the open-source software OpenFOAM v8 was used, see \url{openfoam.org/version/8/}. The C++ Finite-Volume-Method code is a toolbox containing numerical solvers to solve not only but mostly computational fluid dynamics problems. The framework for Parareal was written in Python with the PyFoam, see \url{github.com/dicehub/PyFoam}, and NumPy \cite{NumPy} libraries. For the creation of the meshes the \verb|blockMesh| utility by OpenFoam has been used. The utility requires the dictionary \verb|BlockMeshDict|, in which the computational domain is described by blocks that are defined by Cartesian coordinates, boundary conditions and amount of finite volume cells within the respective blocks. Since the domain in this test case is relatively simple, we decided to write a C++ tool that would generate the \verb|blockMeshDict| allowing for parametrization of the mesh resolution and geometry. 

\subsection{Flow around a cylinder}

For our numerical experiments we chose to simulate laminar flow past a cylinder. The Reynolds number for this test case reaches up to $Re=1,000$ to cover several flow regimes. At $Re=100$ laminar flow with periodic vortex shedding develops, between $180<Re<400$ three-dimensional instabilities cause a streamwise vortex structure and at $Re=1,000$ the transition regime to turbulent flow is approached, compare \cite{Rajani2008}. Although, a laminar flow problem is simulated, we were able to test the Micor-Macro Parareal algorithm against different physical phenomenons by choosing this one test case. For the sake of simplicity in this first assessment of the algorithm we started with a 2-D representation of the cylinder flow. Since there are differences between two and three dimensions in space in the simulation results \cite{So2005}, and therefore in agreement with experimental data sets, we did not aim for best possible correctness of the results. The major goal of this study represents the feasibility of approximating the serial fine simulation by Parareal, as it is introduced in the following. 

In this test case the incompressible Navier-Stokes equations for laminar flow are solved. OpenFOAM provides the transient laminar solver \verb|icoFoam| that solves the continuity equation:
\begin{equation}
	\nabla \cdot \mathbf{u} \; = \; 0 \; , \quad \mathbf{u} \; = \; \left( U_x \; , \; U_y \; , \; U_z \right) \; ,
\end{equation}
and the momentum equation:
\begin{equation}
	\frac{\partial}{\partial t} \mathbf{u} \; + \; \left( \mathbf{u} \cdot \nabla \right) \mathbf{u} \; = \; - \nabla p \; + \; \nu \Delta \mathbf{u} \; ,
\end{equation}
with $\mathbf{u}$ denoting the velocity and $p$ the kinematic pressure, the three velocity components $\left( U_x \; , \; U_y \; , \; U_z \right)$ are given for completeness, as the test case is set up in two dimensions. \verb|icoFoam| uses the PISO algorithm to solve the Navier-Stokes equations. As numerical time integration method for both propagators the backward Euler scheme was chosen. The lift coefficient $C_L$ is defined as:
\begin{equation}
	C_L \; = \; \frac{2 f_L}{\rho \; D \; U_\infty^2 \; A_{ref}} \; 
\end{equation}
where $f_L$ denotes the lift force, $D=2$ the cylinders diameter and $A_{Ref}$ the references surface of the cylinder. The lifting force $f_L$ and lift coefficient $C_L$ is computed by the \verb|forceCoeffs| function provided by OpenFOAM and written to an output file during run-time. 
\begin{figure}[H]
	\centering
	\includegraphics{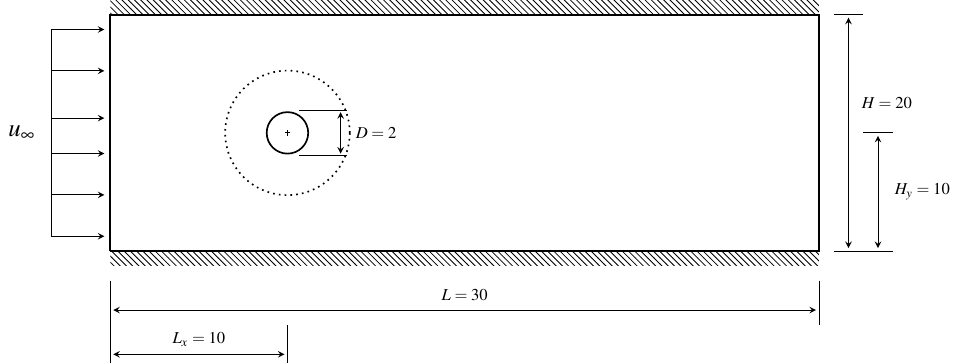}
	\caption{The geometry for the channel flow around a cylinder. $H$ and $L$ denote the height and length of the channel, while $L_x$ and $H_y$ represent the cylinder coordinates in the domain, respectively. The cylinder radius is defined by $R=D/2=1m$. For all numerical experiments the homogeneous inlet velocity is given by $u_\infty=1m/s$. A refinement section was used to increase the spatial resolution around the cylinder (dotted circle).}
	\label{FIG:geometry}
\end{figure}

\begin{figure}[H]
	\centering
	\includegraphics[width=.47\textwidth]{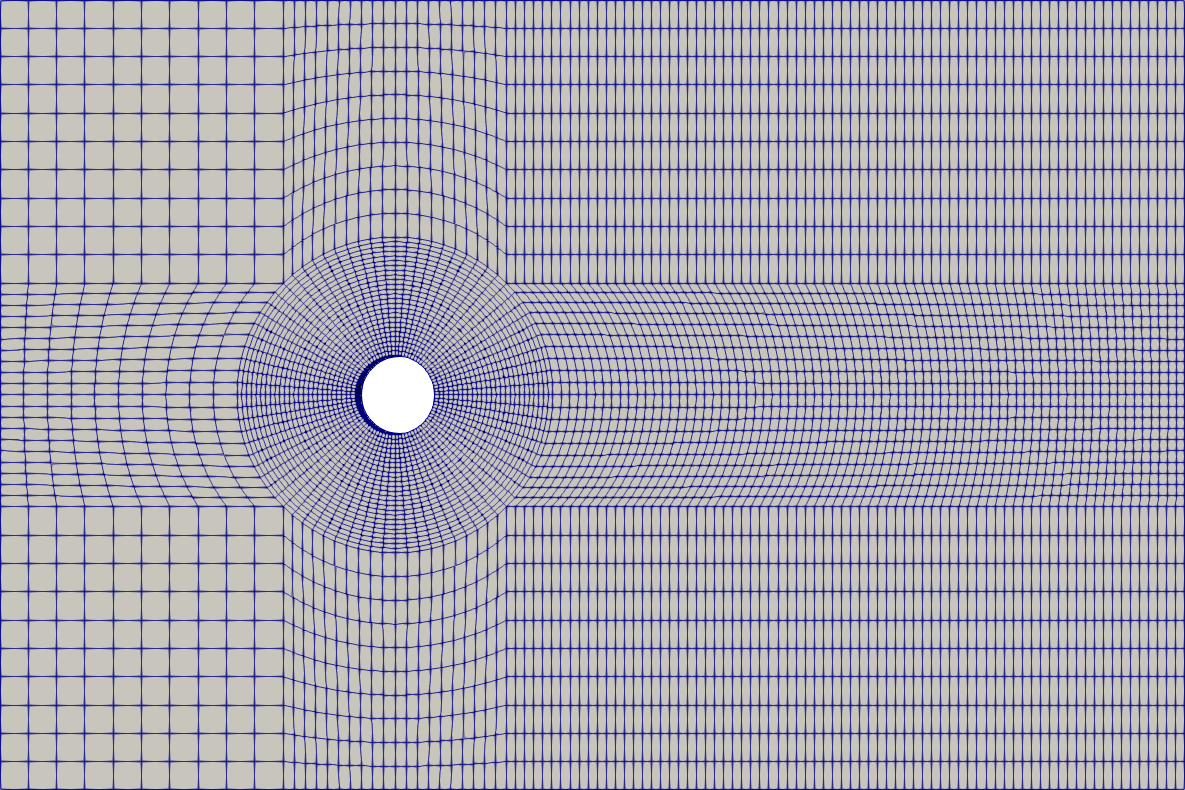}\hfill
	\includegraphics[width=.47\textwidth]{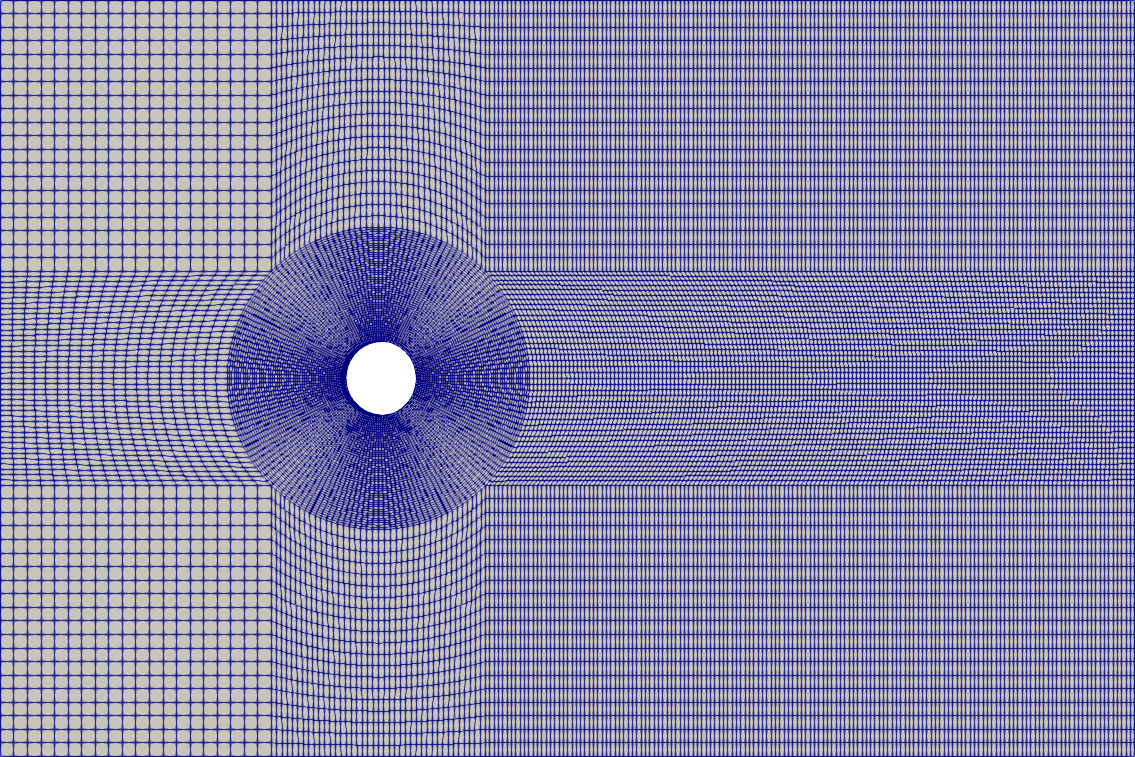}
	\caption{Comparison of the low (left) and high (right) resolution mesh. The coarse mesh contains 5800 cells and the fine mesh 23200 cells.}
	\label{FIG:resolution}
\end{figure}
OpenFOAM provides the preprocessor \verb|blockMesh| for mesh generation. The geometry of the domain is shown in Fig.\ref{FIG:geometry}. We assumed an homogeneous inlet velocity $u_\infty = 1 m/s$. Boundary conditions for cylinder walls are non-slip. On the top and bottom boundary of the domain the symmetric plane boundary condition is imposed. In Fig.\ref{FIG:resolution} the meshes generated by \verb|blockMesh| are depicted. We designed the high resolution mesh to hold four times as many cells as the coarse mesh. The area around is refined in order to obtain better results for the lift coefficient while saving computational effort in the areas in front (stream-wise) of the cylinder. In the created C++ tool for the generation of the mesh, the parameters for the geometry dimensions, amount of finite volumes and refinement zone are freely modifiable by a parameter file.

\subsection{Interpolation Methods}

The \verb|OpenFOAM| software provides three interpolation schemes that can be applied by the \verb|mapFields| utility. \verb|mapNearest| (NN), \verb|interpolate| (IN) and \verb|cellPointInterpolate| (CP) interpolation. Each method is called via command line within the desired destination folder by:
\begin{lstlisting}
	mapFields $SOURCE -sourceTime $TIME -consistent -mapMethod $METHOD
\end{lstlisting}
where \verb|$SOURCE| defines the location of the fields of the source mesh. The source time specifies the solution within the source folder to a given point in time. If the geometry and boundary conditions of target and source mesh is identical, the option \verb|consistent| is set. \verb|mapMethod| selects one of the provided interpolation types, the default is \verb|interpolate|.  The \verb|mapNearest| option corresponds to the nearest neighbor method and  assumes the value of an interpolated cell point by choosing the source mesh cell points value, which is closest in space. With \verb|interpolate| a (bi-/tri-)linear interpolation method is applied. \verb|cellPointInterpolate| represents the inverse volume weighted interpolation method. 

Unfortunately, the \verb|mapFields| utility is not able to interpolate variables that are not stored on the cellpoints. The fluxes are located at the cell faces and therefore are not considered during the mapping of fields between to meshes. If the simulation is continued, the missing fluxes will be generated from the velocity components. This approach violates the continuity and becomes apparent during the evaluation of the lift coefficient in Sec.\ref{SEC:EXP}.

\begin{table}[H]
	\centering
	\caption{Interpolation error estimates for $Re=100$ at $T=200$ seconds.}
	\begin{tabular}{|c|c|c|c|c|}
		\hline
		Interpolation method & Var. & $G - \mathbf{R} \circ \mathbf{L} (G)$ & $F - \mathbf{L} \circ \mathbf{R} (F)$ & $G - \mathbf{R} \circ \mathbf{P}(G)$\\
		\hline 
		& $U_x$ & 0.0 & 4.1e-01 & 0.0 \\
		\verb|mapNearest| & $U_y$ & 0.0 & 2.9e-01 & 0.0 \\
		& $p$ & 0.0 & 8.6e-02 & 0.0 \\
		\hline
		& $U_x$ & 1.3e-01 & 3.8e-01 & 5.9e-02 \\
		\verb|interpolate| & $U_y$ & 6.5e-02 & 2.5e-01 & 3.3e-02 \\
		& $p$ & 3.1e-02 & 6.3e-02 & 2.2e-02 \\
		\hline
		& $U_x$ & 8.5e-02 & 7.0e-02 & 2.7e-02 \\
		\verb|cellPointInterpolate| & $U_y$ & 9.1e-02 & 7.3e-02 & 2.7e-02 \\
		& $p$ & 1.1e-02 & 1.2e-02 & 5.2e-03 \\
		\hline
	\end{tabular}
	\label{TAB:intpRe100}
\end{table}
The consistency requirements, given in Eq.\ref{EQ:requirements}, are given in Tabs.\ref{TAB:intpRe100} and \ref{TAB:intpRe1000} for $Re=100$ and $Re=1,000$. Of the three methods, only \verb|mapNearest| meets the consistency properties. For $Re=100$ the \verb|cellPointInterpolate| method shows better error estimates up to one order compared to the linear interpolation scheme. With increasing Reynolds number both methods equalize. The cell centers in the coarse and fine mesh are at different locations in the domain, and since both methods interpolate the center values the consistency properties cannot be achieved. Although the consistency conditions are fulfilled by \verb|mapNearest|, it should be noted that for interpolation from coarse to fine grids, piece-wise constant functions emerge in some areas. In addition to the before-mentioned inability to interpolate the fluxes, this scheme introduces additional discontinuity.
\begin{table}[H]
	\centering
	\caption{Interpolation error estimates for $Re=1,000$ at $T=200$ seconds.}
	\begin{tabular}{|c|c|c|c|c|}
		\hline
		Interpolation method & Var. & $G - \mathbf{R} \circ \mathbf{L} (G)$ & $F - \mathbf{L} \circ \mathbf{R} (F)$ & $G - \mathbf{R} \circ \mathbf{P}(G)$\\
		\hline 
		& $U_x$ & 0.0 & 6.5e-01 & 0.0\\
		\verb|mapNearest| & $U_y$ & 0.0 & 4.1e-01 & 0.0 \\
		& $p$ & 0.0 & 1.8e-01 & 0.0 \\
		\hline
		& $U_x$ & 2.2e-01 & 3.9e-01 & 2.6e-01\\
		\verb|interpolate| & $U_y$ & 8.2e-02 & 2.8e-01 & 1.3e-01 \\
		& $p$ & 4.5e-02 & 1.2e-01 & 6.3e-02 \\
		\hline
		& $U_x$ & 3.1e-01 & 4.7e-01 & 2.6e-01 \\
		\verb|cellPointInterpolate| & $U_y$ & 2.5e-01 & 3.1e-01 & 2.8e-01 \\
		& $p$ & 2.9e-02 & 1.8e-02 & 3.8e-02 \\
		\hline
	\end{tabular}
	\label{TAB:intpRe1000}
\end{table}

\section{Numerical Experiments}\label{SEC:EXP}

An overview of the simulations carried out in this study are given in Tab.\ref{TAB:settings}. We choose to the simulation end time $T=200$ seconds, since it allows for the development of a Karmarn vortex street and evaluation of the diagnostic lift coefficient $C_L$. The test cases are separated by the Reynolds number and in the first case of $Re=100$ by the interpolation method. Each case is subdivided by the amount of time slices $N_t \in [5,10,20]$ applied by Parareal. The runtime ratio $m=m_S \times m_T$ is measured by the execution times and considers the spatial and temporal resolution. The contribution by spatial resolution is $m_S=4$ for all cases. The run-time ratio computed from the different time step sizes $m_T$ is given in Tab.\ref{TAB:settings}. We listed the theoretical ratio $m$ and the measured ratio $\tilde{m}$ below. For the time step sizes we choose the largest allowable with respect to the CFL condition on both meshes. With increasing Reynolds number the coarse time step size $\Delta t_C$ gets gradually closer to $\Delta t_F$. In case of $Re=1,000$ the ratio $m$ is to small for wall-time reduction by time-parallelism and hence, the micro-macro approach is required.
The Parareal simulations were carried out until half of the maximal iterations $K=N_t/2$ were reached. According to Eq.\ref{EQ:speedup} this would allow for a minimal theoretical speed-up of $S\approx2$.  Exception to this presents the test cases with $N_t=5$, where iterations up to the second last iteration were computed to check for further error reduction.  The measured run-time ratios $\tilde{m}$ demonstrates the good scalability properties of OpenFOAM when compared to the theoretical measure $m$. The measurements denote an average over the wall-times in each experiment. Since the goal of the study was to check the feasibility of a time-parallel algorithm, every simulation was carried out on a single thread. 

\begin{table}[H]
	\caption{Test case settings for the Parareal simulations.}
	\centering
	\begin{tabular}{|c|c|c|c|c|c|c|c|c|}
		\hline
		Re &  $N_t$ & $\tilde{m}$ & $m$ & $m_T$ & $\Delta t_F$ & $\Delta t_C$ & Interpolation method & $T$ \\
		\hline
		100 & & & 8 & 2 & 0.05 & 0.1 & NN &\\
		100 & [5,10,20] & 7.6 & 8 & 2 & 0.05 & 0.1 & CP & 200 \\
		100 & & & 8 & 2 & 0.05 & 0.1 & IN &  \\
		\hline
		200 & [5,10,20] & 9.91 & 10 & 2.5 & 0.02 & 0.05 & NN & 200 \\
		\hline
		400 & [5,10,20] & 9.96 & 10 & 2.5& 0.02 & 0.05 & NN & 200 \\
		\hline
		1000 & [5,10,20] & 4.95 & 5 & 1.25 & 0.02 & 0.025 & NN & 200 \\
		\hline
	\end{tabular}
	\label{TAB:settings}
\end{table}

\section{Results}

\begin{figure}[H]
	\centering
	\includegraphics{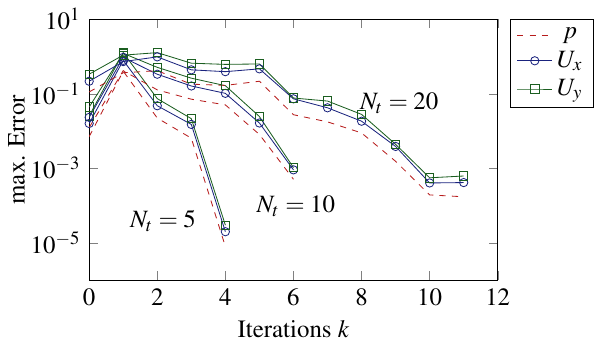}
	\caption{Convergence behavior for $p$, $U_x$ and $U_y$ with respect to the amount of time slices $N_t$. The overall simulation time was $T=200$ seconds. The maximum error estimates are evaluated at the end of the simulation $T=200$. Reynolds number is $Re=100$.}
	\label{FIG:Re100prog}
\end{figure}

The convergence behavior of the prognostic variables with respect to the amount of time slices $N_t$ for the test case $Re=100$ is given in Fig.\ref{FIG:Re100prog}. With increasing $N_t$ the stagnation in error reduction in the first iterations lengthens before the algorithm starts converging. Nevertheless, full convergence to machine precision with respect to the serial reference solution cannot be achieved due to the restart mechanism during the execution of Parareal and the inability to interpolate fluxes. 

\begin{figure}[H]
	\centering
	\begin{subfigure}[t]{.48\textwidth}
		\includegraphics{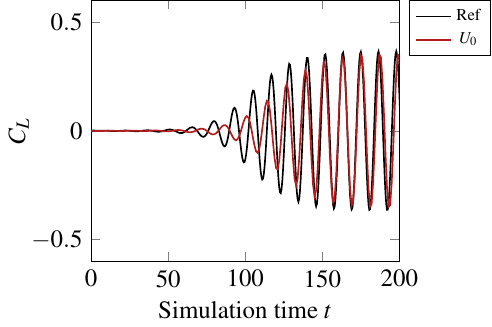}
	\end{subfigure}
	\begin{subfigure}[t]{.48\textwidth}
		\includegraphics{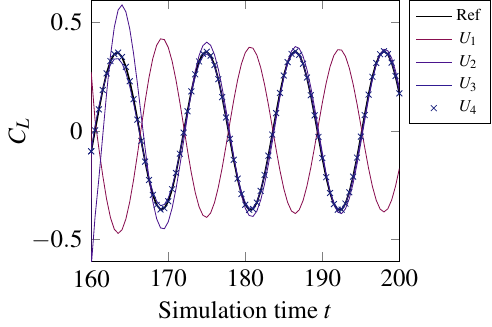}
	\end{subfigure}
	\caption{Left: Lift coefficient $C_L$ of the coarse (red) and fine (black) serial runs for $Re=100$. Right: Parareal convergence of the lift coefficient for the $N_t=5$ case shown for the last time slice.}
	\label{FIG:Re100diag}
\end{figure}

\begin{figure}[H]
	\centering
	\begin{subfigure}[t]{.48\textwidth}
		\includegraphics{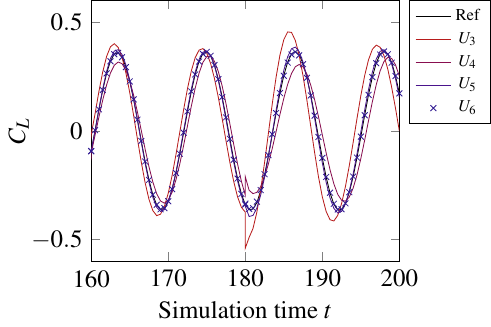}
	\end{subfigure}
	\begin{subfigure}[t]{.48\textwidth}
		\includegraphics{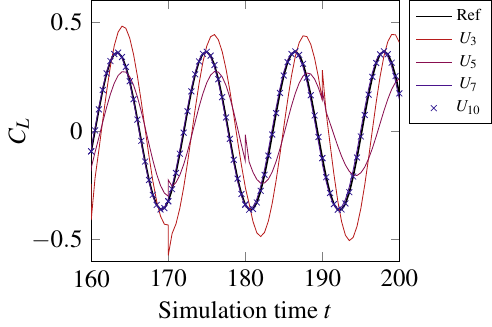}
	\end{subfigure}
	\caption{Parareal convergence of the lift coefficient for the $N_t=10$ (left) and $N_t=20$ (right) cases shown for the respective last time slice. Reynolds number is $Re=100$.}
	\label{FIG:Re100diag2}
\end{figure}

In left panel of Fig.\ref{FIG:Re100diag} the evolution of the lift coefficient $C_L$ is given for the coarse and fine mesh. The lower resolution leads to a delayed and less pronounced development of the Karman vortex street. In right panel the iterations $k=1,2,3$ for the diagnostic lift coefficient $C_L$ are shown. 
Fig.\ref{FIG:Re100diag2} displays the Parareal iterations to the lift coefficient for $N_t=10$ and $N_t=20$. The visible jumps in the progression are discontinuities caused by interpolation and the omission of the fluxes. Although, the convergence of the prognostic variables given in Fig.\ref{FIG:Re100prog} is rather weak, we find the algorithm to converge better with increasing amount of time slices $N_t$. For the last sub-case $N_t=20$ the maximum error reaches $10^{-3}$ in 10 iterations, while acceptable convergence of the diagnostics is reached in 7 iterations.

\begin{figure}[H]
	\centering
	\includegraphics{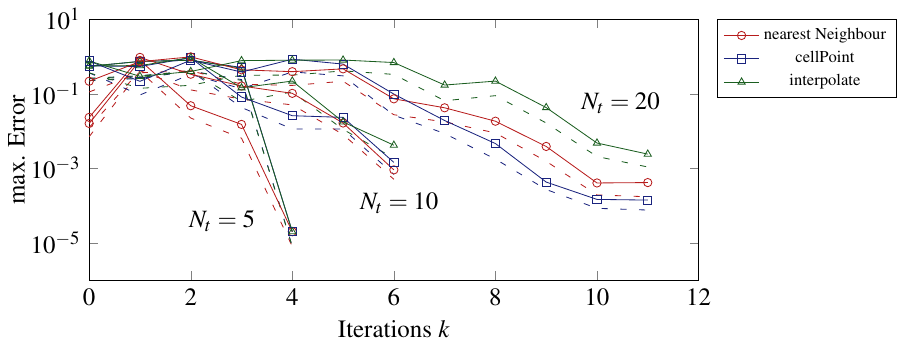}
	\caption{Comparison of interpolation methods with respect to $N_t$ and the prognostic variables $p$ (dashed) and velocity component $U_x$ (solid). Reynolds number is $Re=100$.}
	\label{FIG:Re100intp}
\end{figure}

After the completion of the first test case, where the nearest neighbor scheme was applied, the impact of the other two methods was investigated. A comparison of the convergence to the prognostic variables with respect to the interpolation method is given in Fig.\ref{FIG:Re100intp}. We restricted the visualization to pressure $p$ and the velocity component $U_x$ for the sake of clarity. For the $N_t=5$ sub-case we find the nearest neighbor scheme to be performing best, although all interpolation methods equalize with iteration $k=4$. With increasing $N_t$ the deviations become more distinctive, foremost the linear interpolation scheme demonstrates to be a poor choice for the application in this algorithm. \\

\begin{figure}[H]
	\centering
	\begin{subfigure}[t]{.32\textwidth}
		\includegraphics{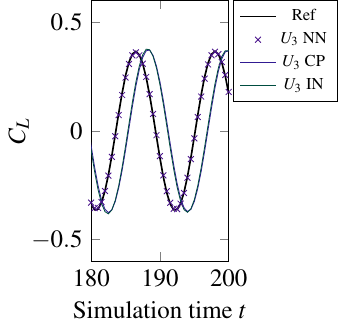}
	\end{subfigure}
	\begin{subfigure}[t]{.32\textwidth}
		\includegraphics{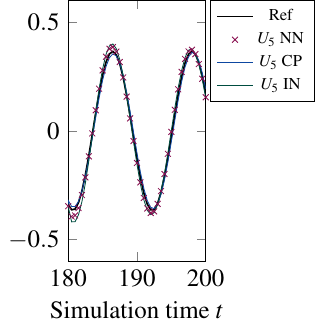}
	\end{subfigure}
	\begin{subfigure}[t]{.32\textwidth}
		\includegraphics{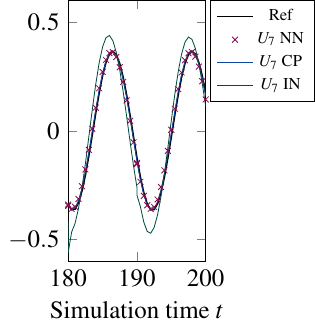}
	\end{subfigure}
	\caption{Left: Parareal convergence of the lift coefficient with respect to the interpolation methods for $N_t=5$. Right: Influence of the interpolation scheme for $N_t=10$. For both cases the iteration before convergence has been chosen. Reynolds number is $Re=100$.}
	\label{FIG:Re100intpDiag}
\end{figure}

In Fig.\ref{FIG:Re100intpDiag} the time course of the lift coefficient is shown for $t \in [180,200]$ seconds. We choose the iterations before the reference result was reproduced in order to point out the difference in the mappings between the meshes. In the left panel the \verb|interpolate| and \verb|cellPointInterpolate| methods we find the time history of the lift coefficient is shifted in time. The nearest neighbor scheme produces, in agreement with the convergence in Fig.\ref{FIG:Re100intp}, the best results for all cases. Thus, the difference between NN and CP diminishes with increasing amount of time slices. The IN method shows overall the poorest performance, again in agreement with the evaluation of the prognostic variables.

\begin{figure}[H]
	\centering
	\includegraphics{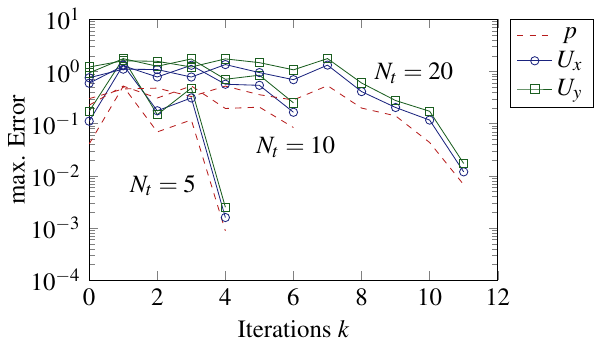}
	\caption{Parareal convergence of prognostic variables with respect to the amount of time slice $N_t$. The error were estimated at the end of the simulation $T=200$ seconds. Reynolds number is $Re=200$.}
	\label{FIG:Re200prog}
\end{figure}

\begin{figure}[H]
	\centering
	\begin{subfigure}[t]{.32\textwidth}
		\includegraphics{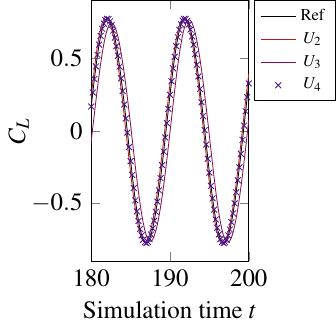}
	\end{subfigure}
	\begin{subfigure}[t]{.32\textwidth}
		\includegraphics{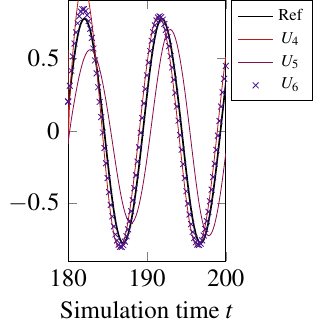}
	\end{subfigure}
	\begin{subfigure}[t]{.32\textwidth}
		\includegraphics{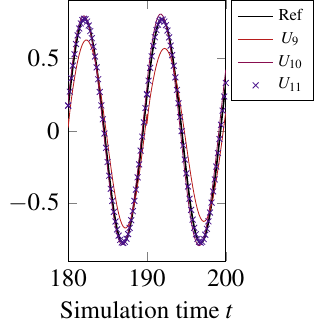}
	\end{subfigure}
	\caption{Convergence of the diagnostic lift coefficient $C_L$ with respect to the amount of time slices $N_t=5$ (left), $N_t=10$ (middle) and $N_t=20$ (right). Reynolds number is $Re=200$. }
	\label{FIG:Re200diag}
\end{figure}

The convergence evaluation of the lift coefficient is given in Fig.\ref{FIG:Re200diag}. In agreement with the poor convergence behavior in Fig.\ref{FIG:Re200prog}, the sub-case $N_t=10$ fails to reproduce the reference for $C_L$. Further iterations would be necessary, which would prohibit the possibility of speed-ups. The other two cases need $k=4$ and $k=11$ iterations, respectively, to converge to the time course reference of the lift coefficient.

\begin{figure}[H]
	\centering
	\includegraphics{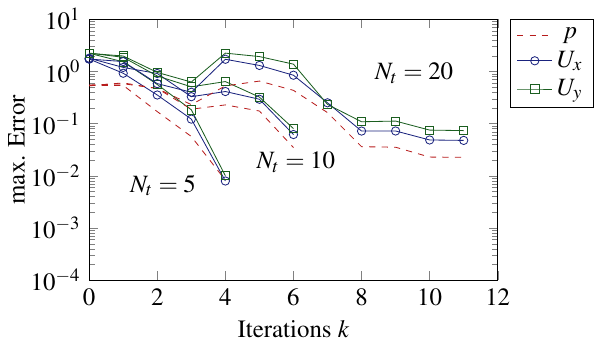}
	\caption{Parareal convergence of prognostic variables with respect to the amount of time slice $N_t$. The error were estimated at the end of the simulation $T=200$ seconds. Reynolds number is $Re=400$.}
	\label{FIG:Re400prog}
\end{figure}

In Fig.\ref{FIG:Re400prog} the convergence results for the case $Re=400$ are given. With further increase of the Reynolds number the convergence error deteriorates. For the last sub-case $N_t=20$ the formation of a bulk is observed between the iterations $k=3$ to $k=6$. The error convergence stagnates around $10^{-1}$ for all sub-cases, and therefore no longer allow the prognostic variables to be considered as sufficiently approximated. 

\begin{figure}[H]
	\centering
	\begin{subfigure}[t]{.32\textwidth}
		\includegraphics{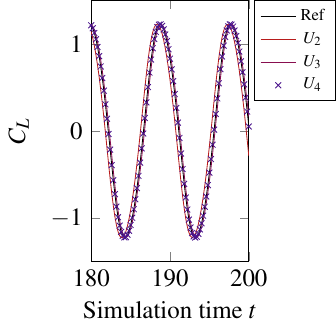}
	\end{subfigure}
	\begin{subfigure}[t]{.32\textwidth}
		\includegraphics{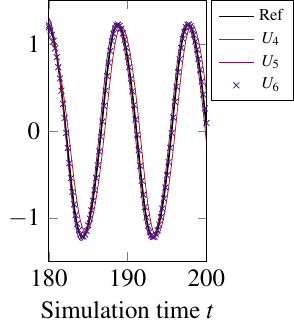}
	\end{subfigure}
	\begin{subfigure}[t]{.32\textwidth}
		\includegraphics{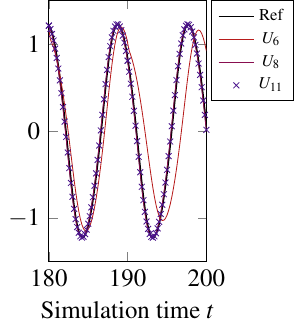}
	\end{subfigure}
	\caption{Convergence of the diagnostic lift coefficient $C_L$ with respect to the amount of time slices $N_t=5$ (left), $N_t=10$ (middle) and $N_t=20$ (right). Reynolds number is $Re=400$.}
	\label{FIG:Re400diag}
\end{figure}

Interestingly, we find the convergence to the diagnostic variable slightly improved in Fig.\ref{FIG:Re400diag}. In the sub-case $N_t=10$ the reference solution can be reproduced within 6 iterations. With $N_t=20$ the reference lift coefficient is sufficiently approximated with $k=8$ iterations. In the case before ($Re=200$) $k=11$ iterations were required. The first sub-case with $N_t=5$ suggests, that convergence in the diagnostic variable can only be achieved, when the algorithm computes more iterations than predefined by the theoretical speed-up estimate in Eq.\ref{EQ:speedup}. 

\begin{figure}[H]
	\centering
	\includegraphics{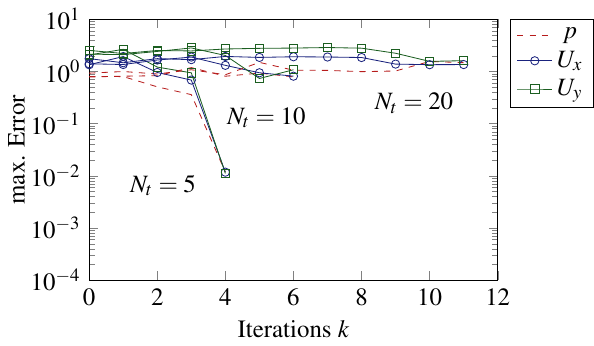}
	\caption{Parareal convergence of prognostic variables with respect to the amount of time slice $N_t$. The error were estimated at the end of the simulation $T=200$ seconds. Reynolds number is $Re=1,000$.}
	\label{FIG:Re1kprog}
\end{figure}

\begin{figure}[H]
	\centering
	\begin{subfigure}[t]{.32\textwidth}
		\includegraphics{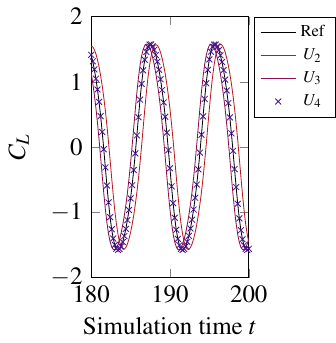}
	\end{subfigure}
	\begin{subfigure}[t]{.32\textwidth}
		\includegraphics{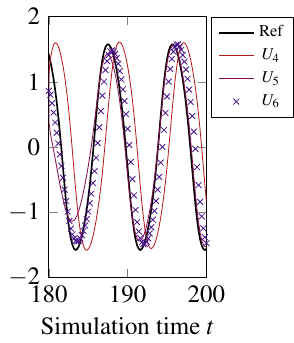}
	\end{subfigure}
	\begin{subfigure}[t]{.32\textwidth}
		\includegraphics{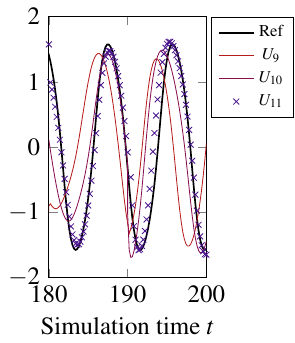}
	\end{subfigure}
	\caption{Convergence of the diagnostic lift coefficient $C_L$ with respect to the amount of time slices $N_t=5$ (left), $N_t=10$ (middle) and $N_t=20$ (right). Reynolds number is $Re=1,000$.}
	\label{FIG:Re1kdiag}
\end{figure}

The evaluation of the prognostic variables of the last test case with $Re=1,000$ is depicted in Fig.\ref{FIG:Re1kprog}. The convergence error stagnates for all iterations in all sub-cases within the chosen amount of maximal iterations $K$. As found in the last numerical experiment, only if the iteration count exceeds the maximal allowable amount of iterations, as demonstrated for $N_t=5$ in the left panel of Fig.\ref{FIG:Re1kdiag}, it is possible to recover the reference solution.  

In agreement with the evaluation of the prognostic variables in Fig.\ref{FIG:Re1kprog}, Parareal failed to converge to the reference time course of the lift coefficient, as given in Fig.\ref{FIG:Re1kdiag}. For all iterations errors in phase and amplitude are still existent. 
With the numerical results given we conclude this section with the following findings:
\vspace{.25cm}
\begin{itemize}
	\item[1.] The convergence deteriorates with increasing Reynolds number. For low Reynolds numbers ($Re\le400$) the algorithm could generate speed-ups, under the condition that the diagnostic variable represents the main objective of the simulation. 
	\item[2.] The algorithm performs better for diagnostic lift coefficient than for the prognostic variables local in space and time. Utilizing more time slices results in better convergence behavior in the diagnostic variable.
	\item[3.] The \verb|mapNearest| scheme proved to be the best choice among the methods provided by OpenFOAM. The requirements regarding the interpolation operators in Eq.\ref{EQ:requirements} of the Micro-Macro Parareal algorithm suggested by \cite{Samaey2012} are fulfilled only in this case. 
	\item[4.] The lack of possibilities to interpolate the fluxes located at the internal cell faces causes discontinuities in the approximation of the lift coefficient and affects the convergence behavior of the algorithm. 
\end{itemize}
\vspace{.25cm}
We further tried to apply Parareal to a turbulent flow regime of $Re=10,000$ by adapting the algorithm to the RANS k-$\varepsilon$ turbulence model. We found the \verb|mapNearest| scheme representing the only possibility for applying Parareal underlining the necessity of the consistency property for interpolations. Nevertheless, the algorithm was not able to compute approximations to a reference solution over the simulation time interval of $T=200$ seconds. Micro-Macro Parareal would suffer from blow-ups in the turbulence kinetic energy $k$. Up to a simulation time of $T=150$ seconds the algorithm could be applied, but the vortex street was not fully developed at that point. Reduction of the time interval and repeated application of Parareal lead to no success either. Hence, we refrained from reporting the results. We conclude from that unsuccessful experiment, that the consistency properties in Eq.\ref{EQ:requirements} of the interpolation operators are in fact important for PDE problems when applying the proposed algorithm. On the other hand provides the \verb|mapNearest| method rather poor approximation quality by adding discontinuities, so that a new interpolation scheme should be implemented that takes into account the operators consistency and conservation of physical quantities.

\section{Conclusion}

In this the micro-macro Parareal algorithm by \cite{Samaey2012}\cite{Slawig2018} was applied to the laminar cylinder flow problem. The time-parallel approximation of the diagnostic lift coefficient was successful for cases of $Re \le 400$. When approaching the transition regime towards turbulent flow the algorithm fails to converge within a reasonable amount of iterations $K\le N_t/2$.  \\
For future work we suggest to address the mapping methods between meshes and the possibility of exploiting the concept of domain decomposition: The interpolation methods available in this work were not able to map the fluxes stored on the internal faces fo the domain. In order to avoid the resulting incontinuities we suggest to implement an interpolation scheme that is capable of mapping the fluxes in a conservative way. Challenging demands are placed on such an implementation, as it must also manage the interpolation between decomposed spatial domains. With such an interpolation method the aspect of parallelism in space and time could be assessed conclusively with respect to the effective wall-time reduction of Parareal. Once these steps have been taken, we expect that the algorithm can be successfully applied to all the presented test cases and for future work to a 3-dimensional realization of the cylinder flow. With an appropriate conservative interpolation scheme the adaption to turbulent regimes with turbulence modeling would offer potential for future investigations. \\
The experiments demonstrated that the consistency property in Eq.\ref{EQ:requirements} for the interpolation methods is important for the convergence behavior of the Micro-Macro Parareal algorithm, though it originally was formulated for ODEs. 


\printbibliography	

\end{document}